\documentclass[11pt,frenchb,english]{article}
\usepackage{palatino}
\usepackage[T1]{fontenc}
\usepackage[latin1]{inputenc}
\usepackage{a4}
\usepackage{amssymb}

\makeatletter

\newcommand{\boldsymbol}[1]{\mbox{\boldmath $#1$}}

\providecommand{\tabularnewline}{\\}

 \usepackage{amsthm, amsmath, amsfonts, amssymb}
 \theoremstyle{plain}    
 \newtheorem{thm}{Theorem}%[section] enlevé par moi pour eviter la numérotation par section
 \theoremstyle{plain}    
 \newtheorem{prop}[thm]{Proposition} %%Delete [thm] to re-start numbering

 \theoremstyle{definition}
  \newtheorem*{example*}{Example}
 \theoremstyle{plain}    
 \newtheorem{lem}[thm]{Lemma} %%Delete [thm] to re-start numbering 
 \theoremstyle{plain}    
 \newtheorem{cor}[thm]{Corollary} %%Delete [thm] to re-start numbering

\usepackage{babel}
\makeatother
\begin{document}

\title{The geometry of nondegeneracy conditions in completely integrable
systems}

\author{Nicolas Roy%
\footnote{\noindent Address : Geometric Analysis Group, Institut für Mathematik,
Humboldt Universität, Rudower Chaussee 25, Berlin D-12489, Germany.
Email : \texttt{roy@math.hu-berlin.de}%
}}

\maketitle
\begin{abstract}
Nondegeneracy conditions need to be imposed in K.A.M. theorems to
insure that the set of diophantine tori has a large measure. Although
they are usually expressed in action coordinates, it is possible to
give a geometrical formulation using the notion of regular completely
integrable systems defined by a fibration of a symplectic manifold
by lagrangian tori together with a Hamiltonian function constant on
the fibers. In this paper, we give a geometrical definition of different
nondegeneracy conditions, we show the implication relations that exist
between them, and we show the uniqueness of the fibration for non-degenerate
Hamiltonians.
\end{abstract}
\textbf{}

\selectlanguage{frenchb}
\begin{abstract}
Dans les théorèmes de type K.A.M., on doit imposer des conditions
de non-dégénérescence pour assurer que l'ensemble des tores diophantiens
a une grande mesure. Elles sont habituellement présentées en coordonnées
actions, mais il est possible d'en donner une formulation géométrique
en considérant des systèmes complètement intégrables définis par la
donnée d'une fibration d'une variété symplectique par des tores lagrangiens
et d'un Hamiltonien constant sur les fibres. Dans cet article, nous
donnons une définition géométrique de différentes conditions de non-dégénérescence,
nous montrons les différentes relations d'implication qui existent
entre elles, et nous montrons l'unicité de la fibration pour les Hamiltoniens
non-dégénérés.\newpage 
\end{abstract}

\selectlanguage{english}
\section*{Introduction}

On a symplectic manifold $\left(\mathcal{M},\omega\right)$, the completely
integrable systems (CI in short) are the dynamical systems defined
by a Hamiltonian $H\in C^{\infty}\left(\mathcal{M}\right)$ admitting
a \emph{momentum map}, i.e. a set $\boldsymbol{A}=\left(A_{1},...,A_{d}\right):\mathcal{M}\rightarrow\mathbb{R}^{d}$
of smooth functions, $d$ being half of the dimension of $\mathcal{M}$,
satisfying $\left\{ A_{j},H\right\} =0$ and $\left\{ A_{j},A_{k}\right\} =0$
for all $j,k:1...d$, and whose differentials $dA_{j}$ are linearly
independent almost everywhere. Then, the Arnol'd-Mineur-Liouville
Theorem \cite{arnold_1,mineur,liouville} insures that in a neighbourhood
of any connected component of any compact regular fiber $\boldsymbol{A}^{-1}\left(a\right)$,
$a\in\mathbb{R}^{d}$, there exists a fibration in lagrangian tori
along which $H$ is constant. These tori are thus invariant by the
dynamics generated by the associated Hamiltonian vector field $X_{H}$. 

Despite the {}``local'' character of the Arnol'd-Mineur-Liouville
Theorem, one might be tempted to try to glue together these {}``local''
fibrations in the case of \emph{regular} Hamiltonians, i.e. those
for which there exists, near each point of \emph{}$\mathcal{M}$,
a local fibration in invariant lagrangian tori. Nevertheless, we would
like to stress the fact that not all regular completely integrable
Hamiltonians are constant along the fibers of a fibration in lagrangian
tori. For example, a free particle moving on the sphere $S^{2}$ can
be described by a Hamiltonian system on the symplectic manifold $T^{*}S^{2}$.
If we restrict ourself to the symplectic manifold $\mathcal{M}=T^{*}S^{2}\setminus S^{2}$,
we can easily show that $\mathcal{M}$ is diffeomorphic to $SO\left(3\right)\times\mathbb{R}$
and that the Hamiltonian $H$ depends only on the second factor. The
energy levels $H=cst$ are thus diffeomorphic to $SO\left(3\right)$.
On the other hand, if there exists a fibration in lagrangian tori
such that $H$ is constant along the fibers, then each energy level
is itself fibered by tori. But a simple homotopy group argument shows
that there exists no fibration of $SO\left(3\right)$ by tori. 

This example actually belongs to the non-generic (within the class
of regular CI Hamiltonians) class of \emph{degenerate} Hamiltonians.
Those Hamiltonians might not admit any (global) fibration in lagrangian
tori, or they might admit several different ones. But, as we will
see, imposing a nondegeneracy condition insures that there exists
a global a fibration of $\mathcal{M}$ in lagrangian tori along which
$H$ is constant, and moreover that it is unique. 

On the other hand, nondegeneracy conditions arise in the K.A.M. theory
where one studies the small perturbations $H+\varepsilon K$ of a
given CI Hamiltonian $H$. The K.A.M. Theorem actually deals with
the \emph{regular part} of a completely integrable system and is usually
expressed in angle-action coordinates. This theorem actually gives
two independent statements. The first statement is that the tori on
which $X_{H}$ verifies a certain diophantine relation are only slightly
deformed and not destroyed by the perturbation $\varepsilon K$, provided
$\varepsilon$ is sufficiently small. The second statement is that
the set of these tori has a large mesure whenever $H$ is non-degenerate. 

There exist different K.A.M. theorems based on different nondegeneracy
conditions, such as the earliest ones of Arnol'd \cite{arnold_1}
and Kolmogorov \cite{kolmogorov}, or those introduced later by Bryuno
\cite{bruno} and Rüssmann \cite{russmann}. They are always presented
in action-angle coordinates and this hides somehow their geometrical
content. But they can be expressed in a geometric way if we consider
CI systems defined on a symplectic manifold $\mathcal{M}$ by a fibration
in lagrangian tori $\mathcal{M}\overset{\pi}{\rightarrow}\mathcal{B}$,
where $\mathcal{B}$ is any manifold, together with a Hamiltonian
$H\in C^{\infty}\left(\mathcal{M}\right)$ constant along the fibers
$\pi^{-1}\left(b\right)$, $b\in\mathcal{B}$. Such a Hamiltonian
must have the form $H=F\circ\pi$, with $F\in C^{\infty}\left(\mathcal{B}\right)$,
and all the nondegeneracy conditions express simply in function of
$F$ and of a torsion-free and flat connection which naturally exists
on the base space $\mathcal{B}$ of the fibration.

In the first section, we review the geometric structures associated
with a fibration in lagrangian tori that allow one to define the connection
on the base space. In Section 2, we give several nondegeneracy conditions,
including those mentionned above, expressed both in a geometric way
and in flat coordinates. Then, we show in Section 3 the different
implication relations that exist between these different conditions.
Finally, in the last section, we give some properties of non-degenerate
CI hamiltonians, as for example the uniqueness of the fibration in
lagrangian tori.

\section{Geometric setting\label{sec_Geometric_setting}}

Let $\left(\mathcal{M},\omega\right)$ be a symplectic manifold of
dimension $2d$ and let $\left(H,\mathcal{M}\overset{\pi}{\rightarrow}\mathcal{B}\right)$
be a regular CI system composed of a fibration in lagrangian tori
$\mathcal{M}\overset{\pi}{\rightarrow}\mathcal{B}$ together with
a Hamiltonian $H\in C^{\infty}\left(\mathcal{M}\right)$ constant
along the fibers $\mathcal{M}_{b}=\pi^{-1}\left(b\right)$, $b\in\mathcal{B}$.
Since by definition the fibers are connected, $H$ must be of the
form $H=F\circ\pi$, with $F\in C^{\infty}\left(\mathcal{B}\right)$.
On the other hand, Duistermaat showed in \cite{duistermaat} that
there exists a natural torsion-free and flat connection $\nabla$
on the base space $\mathcal{B}$ of each fibration in lagrangian tori.
It can be seen as follows. 

First of all, a theorem due to Weinstein \cite{weinstein1,weinstein2}
insures that there exists a natural torsion-free and flat connection
on each leaf of a lagrangian foliation. Moreover, whenever this foliation
defines locally a fibration, then the holonomy of the connection must
vanish. Given a fibration in lagrangian tori $\left(H,\mathcal{M}\overset{\pi}{\rightarrow}\mathcal{B}\right)$,
the space $\mathcal{V}_{\nabla}\left(\mathcal{M}_{b}\right)$ of parallel
vector fields on $\mathcal{M}_{b}=\pi^{-1}\left(b\right)$, for each
$b\in\mathcal{B}$, is thus a vector space of dimension $d$, and
the union $\bigcup_{b\in\mathcal{B}}\mathcal{V}_{\nabla}\left(\mathcal{M}_{b}\right)$
is actually endowed with a structure of a smooth vector bundle over
$\mathcal{B}$. 

On the other hand, since each fiber $\mathcal{M}_{b}$ is a standard%
\footnote{Here, {}``standard'' means holonomy-free.%
} affine torus, one can define for each $b\in\mathcal{B}$ the space
$\Lambda_{b}\subset\mathcal{V}_{\nabla}\left(\mathcal{M}_{b}\right)$
of $1$-periodic parallel vector fields on $\mathcal{M}_{b}$, which
is easily shown to be a lattice of $\mathcal{V}_{\nabla}\left(\mathcal{M}_{b}\right)$.
Moreover, one can show that the union $\Lambda=\bigcup_{b\in\mathcal{B}}\Lambda_{b}$
is a smooth lattice subbundle of $\bigcup_{b\in\mathcal{B}}\mathcal{V}_{\nabla}\left(\mathcal{M}_{b}\right)$,
called the \emph{period bundle}. This is in fact the geometrical content
of the Arnol'd-Mineur-Liouville Theorem \cite{arnold_1,mineur,liouville}.
To prove this, one constructs explicitly smooth sections of $\bigcup_{b\in\mathcal{B}}\mathcal{V}_{\nabla}\left(\mathcal{M}_{b}\right)$
which are $1$-periodic, namely Hamiltonian vector fields $X_{\xi\circ\pi}$
whose Hamiltonian is the pullback of a function $\xi\in C^{\infty}\left(\mathcal{M}\right)$
of a special type and called \emph{action}. Now, the symplectic form
on $\mathcal{M}$ provides for each $b$ an isomorphism between $\mathcal{V}_{\nabla}\left(\mathcal{M}_{b}\right)$
and $T_{b}^{*}\mathcal{B}$. The image of the bundle $\Lambda$ by
this isomorphism is then a smooth lattice subbundle $E^{*}$ of $T^{*}\mathcal{B}$,
called the \emph{Action bundle}, and its dual $E$ (called the \emph{Resonance
bundle}) is a smooth lattice subbundle of $T\mathcal{B}$. This lattice
subbundle $E$ provides a way to associate the tangent spaces $T_{b}B$
for neighbouring points $b$. This thus implies the existence of a
natural integer, torsion-free and flat connection $\nabla$ on the
base space $\mathcal{B}$ (as discovered by Duistermaat \cite{duistermaat}).
Actually, \emph{angle-action} coordinates are semi-global canonical
coordinates $\left(x,\xi\right):\pi^{-1}\left(\mathcal{O}\right)\rightarrow\mathbb{T}^{d}\times\mathbb{R}^{d}$,
where $\mathcal{O}$ is an open subset of $\mathcal{B}$, with the
properties that the $x_{j}$'s are flat (with respect to Weinstein's
connection) coordinates on the tori, and the differentials $d\xi_{j}$
are smooth sections of the Action bundle $E^{*}$ (this implies that
the coordinates $\xi_{j}$ are flat with respect to Duistermaat's
connection on $\mathcal{B}$).

In the sequel, the space of parallel vector fields will be denoted
by $\mathcal{V}_{\nabla}\left(\mathcal{B}\right)$ and the space of
parallel $1$-forms by $\Omega_{\nabla}^{1}\left(\mathcal{B}\right)$.
We mention that in general the holonomy of Duistermaat's connection
does not vanish. As a consequence, the spaces $\mathcal{V}_{\nabla}\left(\mathcal{B}\right)$
and $\Omega_{\nabla}^{1}\left(\mathcal{B}\right)$ might be empty.
Nevertheless, when one works locally in a simply connected subset
$\mathcal{O}\subset\mathcal{B}$, the spaces of local parallel sections
$\mathcal{V}_{\nabla}\left(\mathcal{O}\right)$ and $\Omega_{\nabla}^{1}\left(\mathcal{O}\right)$
are $d$-dimensional vector spaces.

\section{Different nondegeneracy conditions}

Let $\left(H,\mathcal{M}\overset{\pi}{\rightarrow}\mathcal{B}\right)$
be a regular CI system composed of a fibration in lagrangian tori
$\mathcal{M}\overset{\pi}{\rightarrow}\mathcal{B}$ together with
a Hamiltonian $H\in C^{\infty}\left(\mathcal{M}\right)$ constant
along the fibers. As mentionned before, $H$ is of the form $H=F\circ\pi$,
with $F\in C^{\infty}\left(\mathcal{B}\right)$. It turns out that
all the nondegeneracy conditions are expressed in terms of the function
$F$ and Duistermaat's connection $\nabla$ that naturally exists
on $\mathcal{B}$. On the other hand, these conditions are local :
$F$ is said to be \emph{non-degenerate} if is non-degenerate at each
$b\in\mathcal{B}$. Moreover, some of these conditions are expressed
in terms of the space of parallel vector fields $\mathcal{V}_{\nabla}\left(\mathcal{B}\right)$,
but the local character of the nondegeneracy conditions means that
one needs actually only the spaces $\mathcal{V}_{\nabla}\left(\mathcal{O}\right)$
of local parallel vector fields in a neighbourhood $\mathcal{O}\subset\mathcal{B}$
of each point $b\in\mathcal{B}$. We will use a slight misuse of language
and say {}``for each $X\in\mathcal{V}_{\nabla}\left(\mathcal{B}\right)$''
instead of {}``for each $b\in\mathcal{B}$, each neighbourhood $\mathcal{O}\subset\mathcal{B}$
of $b$ and each $X\in\mathcal{V}_{\nabla}\left(\mathcal{O}\right)$''.

For our purposes, let's define for each $X\in\mathcal{V}_{\nabla}\left(\mathcal{B}\right)$
the function $\Omega_{X}\in C^{\infty}\left(\mathcal{B}\right)$ by
$\Omega_{X}=dF\left(X\right)$ and the associated \emph{resonance
set} \[
\Sigma_{X}=\left\{ b\in\mathcal{B}\mid\Omega_{X}\left(b\right)=0\right\} .\]
 We will also denote by $\mathcal{K}=\bigcup_{b}\mathcal{K}_{b}$
the integrable distibution of hyperplanes $\mathcal{K}_{b}\subset T_{b}\mathcal{B}$
tangent to the hypersurfaces $F=cst$, i.e. $\mathcal{K}_{b}=\ker\left.dF\right|_{b}.$
The \emph{Hessian} $\nabla\nabla F$, which is a $\left(0,2\right)$-tensor
field on $\mathcal{B}$, will be denoted by $F^{''}$. It is symmetric
since $\nabla$ is torsion-free. The connection $\nabla$ also yields
an identification of the cotangent spaces $T_{b}^{*}\mathcal{B}$
at neighbouring points $b$ and allows us to define the \emph{frequency
map} $\varphi:\mathcal{B}\rightarrow\Omega_{\nabla}^{1}\left(\mathcal{B}\right)$
by $\varphi\left(b\right)=dF_{b}^{\nabla}$, where $dF_{b}^{\nabla}\in\Omega_{\nabla}^{1}\left(\mathcal{B}\right)$
is the parallel $1$-form which coincides with $dF$ at the point
$b$. In the sequel, the expression $A\propto B$ means that the vectors
$A$ and $B$ are linearly dependent.

We now review different nondegeneracy conditions, including those
used in the literature. We give both a geometrical formulation and
the corresponding (usual) formulation in flat coordinates. 

\bigskip{}
\noindent \textbf{Condition {}``Kolmogorov''} : For each $b\in\mathcal{B}$,
the Hessian, seen as a linear map $F_{b}'':T_{b}\mathcal{B}\rightarrow T_{b}^{*}\mathcal{B}$,
is invertible, i.e \[
X\in T_{b}\mathcal{B},\nabla_{X}\nabla F=0\Longrightarrow X=0.\]
 In flat coordinates, this condition reads : \[
\det\left(\frac{\partial^{2}F}{\partial\xi_{j}\partial\xi_{k}}\right)\neq0.\]

\bigskip{}
\noindent \textbf{Condition {}``Locally diffeomorphic frequency map''}
: The frequency map is a local diffeomorphism. In flat coordinates,
this condition means that the map $\varphi:\xi_{j}\in\mathbb{R}^{d}\rightarrow\frac{\partial F}{\partial\xi_{k}}\in\mathbb{R}^{d}$
is a local diffeomorphism.

\bigskip{}
\noindent \textbf{Condition {}``Iso-energetic'' or {}``Arnol'd''}
: For each $b\in\mathcal{B}$ the restriction to $\mathcal{K}_{b}$
(for the two slots) of the Hessian $\left.F''\right|_{\mathcal{K}_{b}}:\mathcal{K}_{b}\rightarrow\mathcal{K}_{b}^{*}$
is invertible, i.e.\[
X\in\mathcal{K}_{b},\nabla_{X}\nabla F\propto dF\Longrightarrow X=0\]
 In flat coordinates, this condition reads :\[
\det\left(\begin{array}{cc}
\left[\begin{array}{ccc}
\ddots\\
 & \frac{\partial^{2}F}{\partial\xi_{j}\partial\xi_{k}}\\
 &  & \ddots\end{array}\right] & \left[\begin{array}{c}
\vdots\\
\frac{\partial F}{\partial\xi_{j}}\\
\vdots\end{array}\right]\\
\left[\begin{array}{ccc}
\cdots & \frac{\partial F}{\partial\xi_{k}} & \cdots\end{array}\right] & 0\end{array}\right)\neq0.\]

\bigskip{}
\noindent \textbf{Condition {}``Bryuno''} : For each $b\in\mathcal{B}$,
the set of vectors $X\in T_{b}\mathcal{B}$ satisfying $\nabla_{X}\nabla F\propto dF$
is $1$-dimensionnal, i.e. : \[
X,Y\in T_{b}\mathcal{B},\nabla_{X}\nabla F\propto dF\textrm{ and }\nabla_{Y}\nabla F\propto dF\Longrightarrow X\propto Y.\]
This amounts to requiring that for each $b$, the linear map $\mathbf{U}:T_{b}\mathcal{B}\oplus\mathbb{R}\rightarrow T_{b}^{*}\mathcal{B}$
defined by \[
\mathbf{U}\left(X,\alpha\right)=\nabla_{X}\nabla F+\alpha dF\]
 has a rank equal to $d$. In flat coordinates, this condition reads
:\[
\textrm{rank}\left(\begin{array}{cc}
\left[\begin{array}{ccc}
\ddots\\
 & \frac{\partial^{2}F}{\partial\xi_{j}\partial\xi_{k}}\\
 &  & \ddots\end{array}\right] & \left[\begin{array}{c}
\vdots\\
\frac{\partial F}{\partial\xi_{j}}\\
\vdots\end{array}\right]\end{array}\right)=d.\]

\bigskip{}
\noindent \textbf{Condition ''N''} : For each $b\in\mathcal{B}$,
the restriction to $\mathcal{K}_{b}$ (for the first slot) of the
Hessian $\left.F''\right|_{\mathcal{K}_{b}}:\mathcal{K}_{b}\rightarrow T_{b}^{*}\mathcal{B}$
is injective, i.e. : \[
X\in\mathcal{K}_{b},\nabla_{X}\nabla F=0\Longrightarrow X=0.\]
 This is equivalent to requiring that for each $b$, the linear map
$\mathbf{V}:T_{b}\mathcal{B}\rightarrow T_{b}^{*}\mathcal{B}\oplus\mathbb{R}$
defined by \[
\mathbf{V}\left(X\right)=\left(\nabla_{X}\nabla F,dF\left(X\right)\right)\]
has a rank equal to $d$. In flat coordinates, this condition reads
:\[
\textrm{rang}\left(\begin{array}{c}
\left[\begin{array}{ccc}
\ddots\\
 & \frac{\partial^{2}F}{\partial\xi_{j}\partial\xi_{k}}\\
 &  & \ddots\end{array}\right]\\
\left[\begin{array}{ccc}
\cdots & \frac{\partial F}{\partial\xi_{k}} & \cdots\end{array}\right]\end{array}\right)=d.\]

\bigskip{}
\noindent \textbf{Condition {}``Turning frequencies''} : Let $P\left(\Omega_{\nabla}^{1}\left(\mathcal{B}\right)\right)$
denote the projective space of $\Omega_{\nabla}^{1}\left(\mathcal{B}\right)$
and $\pi:\Omega_{\nabla}^{1}\left(\mathcal{B}\right)\rightarrow P\left(\Omega_{\nabla}^{1}\left(\mathcal{B}\right)\right)$
the associated projection. We require the map $\pi\circ\varphi:\mathcal{B}\rightarrow P\left(\Omega_{\nabla}^{1}\left(\mathcal{B}\right)\right)$
to be a submersion. In flat coordinates, this condition amounts to
requiring that the map $\pi\circ\varphi:\mathbb{R}^{d}\rightarrow P\left(\mathbb{R}^{d}\right)$
defined by $\xi\rightarrow\left[\frac{\partial F}{\partial\xi^{j}}\left(\xi\right)\right]$
is a submersion.

\bigskip{}
\noindent \textbf{Condition {}``Iso-energetic turning frequencies''}
: The restriction of the frequency map $\varphi:\mathcal{B}\rightarrow\Omega_{\nabla}^{1}\left(\mathcal{B}\right)$
to each energy level $S_{E}=\left\{ b\mid F\left(b\right)=E\right\} $
is a local diffeomorphism $\pi\circ\varphi$ between $S_{E}$ and
$P\left(\Omega_{\nabla}^{1}\left(\mathcal{B}\right)\right)$.

\bigskip{}
\noindent \textbf{Condition {}``Regular resonant set''} : For each
point $b\in\mathcal{B}$ and each non-vanishing parallel vector field
$X\in\mathcal{V}_{\nabla}\left(\mathcal{B}\right)$, one has \[
d\left(\Omega_{X}\right)_{b}\neq0.\]
 This implies that for each $X\in\mathcal{V}_{\nabla}\left(\mathcal{B}\right)$
the resonant set $\Sigma_{X}$ is a $1$-codimensionnal submanifold
of $\mathcal{B}$.

\bigskip{}
\noindent \textbf{Condition {}``Resonant set with empty interior
''} : For each non-vanishing parallel vector field $X\in\mathcal{V}_{\nabla}\left(\mathcal{B}\right)$,
the resonant set $\Sigma_{X}$ has an empty interior.

\bigskip{}
\noindent \textbf{Condition {}``Rüssmann''} : For each non-vanishing
parallel vector field $X\in\mathcal{V}_{\nabla}\left(\mathcal{B}\right)$,
the image of the frequency map does not annihilate $X$ on an open
set. In flat coordinates, this condition means that the image of $\varphi:\xi_{j}\in\mathbb{R}^{d}\rightarrow\frac{\partial F}{\partial\xi_{k}}\in\mathbb{R}^{d}$
does not lie in any hyperplane passing through the origin.

\section{Hierarchy of conditions}

First of all, we will show that 

\begin{center}\begin{tabular}{cccc|c|c|c|}
\cline{1-1} \cline{3-3} \cline{5-5} \cline{7-7} 
\multicolumn{1}{|c|}{Turning frequencies}&
\multicolumn{1}{c|}{$\Leftrightarrow$}&
\multicolumn{1}{c|}{Bryuno}&
\multicolumn{1}{c|}{$\Leftrightarrow$}&
\multicolumn{1}{c|}{N}&
$\Leftrightarrow$&
\multicolumn{1}{c|}{Regular resonant set}\tabularnewline
\cline{1-1} \cline{3-3} \cline{5-5} \cline{7-7} 
\end{tabular}\end{center}

Consequently, those four equivalent conditions will be denoted by
{}``\emph{Weak nondegeneracy}''. Then, we will show the following
equivalences.

\begin{center}\begin{tabular}{ccc}
\cline{1-1} \cline{3-3} 
\multicolumn{1}{|c|}{Iso-energetic}&
\multicolumn{1}{c|}{$\Longleftrightarrow$}&
\multicolumn{1}{c|}{Iso-energetic turning frequencies}\tabularnewline
\cline{1-1} \cline{3-3} 
\end{tabular}\\
\end{center}

\begin{center}\begin{tabular}{ccc}
\cline{1-1} \cline{3-3} 
\multicolumn{1}{|c|}{Kolmogorov}&
\multicolumn{1}{c|}{$\Longleftrightarrow$}&
\multicolumn{1}{c|}{Locally diffeomorphic frequency map}\tabularnewline
\cline{1-1} \cline{3-3} 
\end{tabular}\\
\end{center}

\begin{center}\begin{tabular}{ccc}
\cline{1-1} \cline{3-3} 
\multicolumn{1}{|c|}{Rüssmann}&
\multicolumn{1}{c|}{$\Longleftrightarrow$}&
\multicolumn{1}{c|}{Resonant set with empty interior}\tabularnewline
\cline{1-1} \cline{3-3} 
\end{tabular}\\
\end{center}

In the second subsection, we will show that we have the following
implications :

\begin{center}\begin{tabular}{cccc|c|}
\cline{5-5} 
\begin{tabular}{|c|}
\hline 
Rüssmann\tabularnewline
\hline
\end{tabular}&
$\Longleftarrow$&
\begin{tabular}{|c|}
\hline 
Weak\tabularnewline
\hline
\end{tabular} &
$\Longleftrightarrow$&
\multicolumn{1}{c|}{\begin{tabular}{c}
Kolmogorov\tabularnewline
or\tabularnewline
Iso-energetic\tabularnewline
\end{tabular}}\tabularnewline
\cline{5-5} 
\end{tabular}\end{center}

\subsection{Equivalent conditions}

\begin{prop}
Condition ''Bryuno'' is equivalent to Condition ''N''.
\end{prop}
\begin{proof}
Consider the linear maps $\mathbf{U}:T_{b}\mathcal{B}\oplus\mathbb{R}\rightarrow T_{b}^{*}\mathcal{B}$
and $\mathbf{V}:T_{b}\mathcal{B}\rightarrow T_{b}^{*}\mathcal{B}\oplus\mathbb{R}$
of Conditions ''Bryuno'' and ''N''. We will show that $\mathbf{U}^{t}=\mathbf{V}$.
Indeed, for each $X\in T_{b}\mathcal{B}$ and each $\left(Y,\alpha\right)\in T_{b}\mathcal{B}\oplus\mathbb{R}$,
the transposed $\mathbf{U}^{t}:T_{b}\mathcal{B}\rightarrow T_{b}^{*}\mathcal{B}\oplus\mathbb{R}$
verifies :\begin{eqnarray*}
\mathbf{U}^{t}\left(X\right)\left(Y,\alpha\right) & = & \mathbf{U}\left(Y,\alpha\right)\left(X\right)\\
 & = & F^{''}\left(Y,X\right)+\alpha dF\left(X\right)\\
 & = & F^{''}\left(X,Y\right)+\alpha dF\left(X\right),\end{eqnarray*}
 where the symmetry of $F^{''}$ has been used. Therefore, one has
\[
\mathbf{U}^{t}\left(X\right)\left(Y,\alpha\right)=\mathbf{V}\left(X\right)\left(Y,\alpha\right).\]
This implies that $\textrm{rank}\left(\mathbf{V}\right)=\textrm{rank}\left(\mathbf{U}^{t}\right)=\textrm{rank}\left(\mathbf{U}\right)$
and thus that Conditions ''N'' and ''Bryuno'' are equivalent.
\end{proof}
\begin{prop}
\label{prop_bruno_ft}Condition ''Bryuno'' is equivalent to Condition
''Turning frequencies''.
\end{prop}
\begin{proof}
Let $dF_{b}^{\nabla}=\varphi\left(b\right)$ be the parallel $1$-form
which coincides with $dF$ at the point $b$. Condition {}``TF''
(''Turning frequencies'') means that the derivative $\left(\pi\circ\varphi\right)_{*}$
of the map $\pi\circ\varphi:\mathcal{B}\rightarrow P\left(\Omega_{\nabla}^{1}\left(\mathcal{B}\right)\right)$
is surjective. Since $P\left(\Omega_{\nabla}^{1}\left(\mathcal{B}\right)\right)$
is $\left(d-1\right)$-dimensionnal, Condition {}``TF'' means that
the kernel of $\left(\pi\circ\varphi\right)_{*}$ is $1$-dimensionnal.
Now, the kernel of $\left(\pi\circ\varphi\right)_{*}$ is precisely
the space of $X$ such that $\varphi_{*}\left(X\right)$ is in the
kernel of $\pi_{*}$, i.e. tangent to the fibers $\pi^{-1}$. Using
the natural isomorphism between $\Omega_{\nabla}^{1}\left(\mathcal{B}\right)$
and its tangent space $T\left(\Omega_{\nabla}^{1}\left(\mathcal{B}\right)\right)$,
one easily sees that $\ker\left(\pi\circ\varphi\right)_{*}$ is composed
of the vectors $X$ such that $\varphi_{*}X_{b}\propto\varphi\left(b\right)$.
Condition {}``TF'' is thus equivalent to requiring that if $X_{b}$
and $Y_{b}$ satisfy $\varphi_{*}X_{b}\propto\varphi\left(b\right)$
and $\varphi_{*}Y_{b}\propto\varphi\left(b\right)$, then $X_{b}\propto Y_{b}$.

On the other hand, we show that $\varphi_{*}X_{b}=\nabla_{X_{b}}\nabla F$.
Indeed, let $t\rightarrow b\left(t\right)$ be a geodesic $t\rightarrow b\left(t\right)$,
passing through $b$ at $t=0$, and let $X_{b}$ the tangent vector
of $b\left(t\right)$ at $b$. In a neigbourhood of $b$, one can
extend $X_{b}$ in a unique way to a parallel vector field $X$. We
thus have $\phi_{X}^{t}\left(b\right)=b\left(t\right)$. We want to
calculate $\varphi_{*}X_{b}=\frac{d}{dt}\left(\varphi\left(b\left(t\right)\right)\right)_{t=0}$.
By definition, for each $t$, $\varphi\left(b\left(t\right)\right)=dF_{b\left(t\right)}^{\nabla}$
is the parallel $1$-form which coincides with $dF_{b\left(t\right)}$
at the point $b\left(t\right)$. It is invariant by the flow of any
parallel vector field, and thus by the one of $X$, i.e. \[
\left(\varphi\left(b\left(t\right)\right)\right)_{b}=\left(\phi_{X}^{t}\right)^{*}\left(dF_{b\left(t\right)}^{\nabla}\right)_{b\left(t\right)}=\left(\left(\phi_{X}^{t}\right)^{*}dF\right)_{b}.\]
 Using the definition of the Lie derivative, one obtains $\varphi_{*}X_{b}=\left(\mathcal{L}_{X}\left(dF\right)\right)_{b}$.
The Cartan's magic formula then gives $\varphi_{*}X_{b}=d\left(dF\left(X\right)\right)_{b}$.
Finally, since $X$ is parallel, we find that $\varphi_{*}X_{b}=\left(\nabla_{X}\nabla F\right)_{b}$. 

Together with the previous result, we have thus proved that is $X_{b}$
and $X_{b}$ satisfy $\left(\nabla_{X_{b}}\nabla F\right)_{b}\propto dF_{b}$
et $\left(\nabla_{Y_{b}}\nabla F\right)_{b}\propto dF_{b}$, then
they must be linearly dependent. This is precisely Condition {}``Bryuno''.
\end{proof}
\begin{prop}
Condition ''N'' is equivalent to Condition ''Regular resonant set''.
\end{prop}
\begin{proof}
Indeed, for each $X\in\mathcal{V}_{\nabla}\left(\mathcal{B}\right)$,
one has $\Omega_{X}=\nabla_{X}F$ and thus $d\left(\Omega_{X}\right)=\nabla\nabla_{X}F=\nabla_{X}\nabla F$.
Conditions ''''Regular resonant set'' thus reads \[
\forall X\in\mathcal{V}_{\nabla}\left(\mathcal{B}\right),\forall b\textrm{ with }X\in\mathcal{K}_{b}\Longrightarrow\left(\nabla_{X}\nabla F\right)_{b}\neq0.\]
This is equivalent to \[
\forall b,\forall X\in\mathcal{K}_{b}\Longrightarrow\left(\nabla_{X}\nabla F\right)_{b}\neq0,\]
i.e. precisely Condition ''N''.
\end{proof}
\begin{prop}
Condition ''Iso-energetic'' is equivalent to Condition ''Iso-energetic
turning frequencies''.
\end{prop}
\begin{proof}
Condition ''Iso-energetic turning frequencies'' amounts to requiring
that at each point $b$ the map $\left(\pi\circ\varphi\right)_{*}$
restricted to an energy level $S_{E}=F^{-1}\left(E\right)$ is a isomorphism
between $T_{b}S_{E}$ and $T_{\pi\left(\varphi\left(b\right)\right)}P\left(\Omega_{\nabla}^{1}\left(\mathcal{B}\right)\right)$,
i.e. the kernel of $\left(\pi\circ\varphi\right)_{*}$ is transverse
to $S_{E}$. Arguing as is the proof of Proposition \ref{prop_bruno_ft},
this is equivalent to requiring that if $X_{b}$ verifies $\left(\nabla_{X_{b}}\nabla H\right)_{b}\propto dH_{b}$,
then $X_{b}$ must be transverse to $S_{E}$, i.e. $dF\left(X_{b}\right)\neq0$,
which is precisely Condition {}``Iso-energetic''.
\end{proof}
\begin{prop}
Condition ''Russmann'' is equivalent to Condition ''Resonant set
with empy interior''.
\end{prop}
\begin{proof}
By definition of the frequency map $\varphi$, for each $X\in\mathcal{V}_{\nabla}\left(\mathcal{B}\right)$,
one has $\varphi\left(b\right)\left(X\right)=dF_{b}^{\nabla}\left(X\right)$.
Since both $X$ and $dF_{b}^{\nabla}$ are parallel, the contraction
$dF_{b}^{\nabla}\left(X\right)$ is a constant function on $\mathcal{B}$
and thus equals to its value, e.g. at the point $b$ which is nothing
but $dF_{b}\left(X\right)=\Omega_{X}\left(b\right)$. We thus have
$\varphi\left(.\right)\left(X\right)=\Omega_{X}\left(.\right)$. Now,
$\varphi\left(b\right)$ annihilates $X$ on an open subset iff the
resonant set $\Sigma_{X}=\Omega_{X}^{-1}\left(0\right)$ contains
this open subset, and therefore its interior is not empty. 
\end{proof}

\subsection{Stronger and weaker conditions}

\begin{prop}
Condition ''Iso-energetic'' implies Condition ''Weak''.
\end{prop}
\begin{proof}
Condition ''Iso-energetic'' means that for each $X\in\mathcal{K}_{b}$,
the 1-form $\left.\nabla_{X}\nabla F\right|_{\mathcal{K}_{b}}$ is
non-vanishing. This property remains true whithout the restriction
to $\mathcal{K}_{b}$, i.e. \[
X\in\mathcal{K}_{b}\Longrightarrow\nabla_{X}\nabla F\neq0,\]
which is Condition ''N''.
\end{proof}
\begin{prop}
Condition ''Kolmogorov'' implies Condition ''Weak''.
\end{prop}
\begin{proof}
Condition ''Kolmogorov'' means that for each $X\in T_{b}\mathcal{B}$,
one has $\nabla_{X}\nabla F\neq0$. By restriction, this remains true
for each $X\in\mathcal{K}_{b}$, which is Condition ''N''.
\end{proof}
\begin{prop}
Condition ''Weak'' implies ''Kolmogorov or Iso-energetic''.
\end{prop}
\begin{proof}
We will actually show the following equivalent logical statement :
if Condition {}``Weak'' is fulfilled but not Condition {}``Kolmogorov'',
then Condition {}``Iso-energetic'' is fulfilled. Suppose that there
is a vector $X\in T_{b}\mathcal{B}$ such that $\nabla_{X}\nabla F=0$
at $b$. We can then extend $X$ around $b$ to a parallel vector
field, and thanks to the symmetry of $F^{''}$, one has $\nabla\nabla_{X}F=0$
at $b$, i.e. $\left(d\Omega_{X}\right)_{b}=0$. On the other hand,
if Condition {}``Regular resonant set'' is fulfilled, this implies
that $b$ cannot belong to the resonant set $\Sigma_{X}$, and thus
$X\notin\mathcal{K}_{b}$. Moreover, Condition ''Bryuno'' insures
that each $Y$ satisfying $\nabla_{X}\nabla F\propto dF$ at $b$
must be linearly dependent on $X$, i.e. $Y\propto X$. We thus have
showed that for each $Y\in\mathcal{K}_{b}$ satisfying $\nabla_{Y}\nabla F\propto dF$,
one has $Y\propto X$ and therefore $Y\notin\mathcal{K}_{b}$. This
implies that $Y=0$. It is precisely Condition ''Iso-energetic''. 
\end{proof}
\begin{prop}
Condition ''Weak'' implies Condition ''Rüssmann''.
\end{prop}
\begin{proof}
Indeed, Condition {}``Regular resonant set'' implies that for each
non-vanishing parallel vector field $X\in\mathcal{V}_{\nabla}\left(\mathcal{B}\right)$,
the resonant set $\Sigma_{X}$ is a $1$-codimensional submanifold,
and thus has an empty interior.
\end{proof}

\subsection{Examples }

\begin{example*}
\noindent \textbf{{}``Kolmogorov\char`\"{} and {}``Iso-energetic}''.
On $\mathcal{B}=\mathbb{R}^{d}\setminus0$, consider the function
$F\left(\xi\right)=\frac{1}{2}\left|\xi\right|^{2}$, where $\left|\xi\right|^{2}=\sum_{j=1}^{d}\left(\xi_{j}\right)^{2}$.
The differential is $dF=\sum_{j=1}^{d}\xi_{j}d\xi_{j}$ and the Hessian
$F_{ij}\left(\xi\right)=\delta_{ij}$ is the identity matrix at each
point $\xi$. . Therefore, one has $\det\left(F_{ij}\right)=1$ ,
which means that $F$ actually satisfies Condition {}``Kolmogorov''
on $\mathbb{R}^{d}$. On the other hand, a straightforward calculation
yields\[
\det\left(\begin{array}{cc}
\left[\begin{array}{ccc}
\ddots\\
 & \frac{\partial^{2}F}{\partial\xi_{j}\partial\xi_{k}}\\
 &  & \ddots\end{array}\right] & \left[\begin{array}{c}
\vdots\\
\frac{\partial F}{\partial\xi_{j}}\\
\vdots\end{array}\right]\\
\left[\begin{array}{ccc}
\cdots & \frac{\partial F}{\partial\xi_{k}} & \cdots\end{array}\right] & 0\end{array}\right)=-\left|\xi\right|^{2}.\]
 This is non-zero whenever $\xi\neq0$, and thus $F$ satisfies both
Conditions {}``Kolmogorov'' and ''Iso-energetic'' on $\mathcal{B}=\mathbb{R}^{d}\setminus0$.
\end{example*}
\noindent \textbf{}

\begin{example*}
\noindent \textbf{{}``Iso-energetic'' but not {}``Kolmogorov''.}
On $\mathcal{B}=\mathbb{R}^{d}\setminus0$, let us consider the function
$F\left(\xi\right)=\left|\xi\right|$. The differential is $dF=\frac{\sum_{j=1}^{d}\xi_{j}d\xi_{j}}{\left|\xi\right|}$
and the Hessian is $F_{ij}\left(\xi\right)=\frac{\delta_{ij}}{\left|\xi\right|}-\frac{\xi_{i}\xi_{j}}{\left|\xi\right|^{3}}$.
We can see that Condition ''Kolmogorov'' is nowhere satisfied since
for each $\xi$, the vector $X_{j}=\xi_{j}$ verifies $\nabla_{X}\nabla F=0$.
Indeed, for each $i$ one has \[
\sum_{j}F_{ij}X_{j}=\frac{\xi_{i}}{\left|\xi\right|}-\frac{\xi_{i}\left|\xi\right|^{2}}{\left|\xi\right|^{3}}=0.\]
Nevertheless, Condition ''Iso-Energetic'' is satisfied since whenever
a vector $X$ verifies $\nabla_{X}\nabla F\propto\nabla F$ and $\nabla_{X}F=0$,
this means that \[
\left\{ \begin{array}{l}
X_{i}-\sum_{j}\frac{\xi_{i}\xi_{j}X_{j}}{\left|\xi\right|^{2}}=\lambda\xi_{i}\\
\sum_{j}X_{j}\xi_{j}=0.\end{array}\right.\]
 By inserting the second equation into the first one, one must have\[
\left\{ \begin{array}{l}
X_{i}=\lambda\xi_{i}\\
\sum_{j}X_{j}\xi_{j}=0,\end{array}\right.\]
 and this is possible only for $\xi=0$, i.e. Condition ''Iso-Energetic''
is satisfied on $\mathcal{B}=\mathbb{R}^{d}\setminus0$.
\end{example*}
\noindent 

\begin{example*}
\textbf{{}``Kolmogorov'' but not {}``Iso-energetic''.} On $\mathcal{B}=\mathbb{R}^{+}\times\mathbb{R}$,
consider the function $F\left(\xi\right)=\frac{\xi_{1}^{3}}{3}+\frac{\xi_{2}^{2}}{2}$.
The differential is $dF=\xi_{1}^{2}d\xi_{1}+\xi_{2}d\xi_{2}$ and
the Hessian is $F^{''}\left(\xi\right)=\left(\begin{array}{cc}
2\xi_{1} & 0\\
0 & 1\end{array}\right)$. The determinant of $F^{''}$ is thus simply $\det\left(F^{''}\right)=2\xi_{1}$
which is non-zero on $\mathcal{B}$. Condition {}``Kolmogorov''
is thus satisfied. Nevertheless, one easily verifies that\[
\det\left(\begin{array}{ccc}
2\xi_{1} & 0 & \xi_{1}^{2}\\
0 & 1 & \xi_{2}\\
\xi_{1}^{2} & \xi_{2} & 0\end{array}\right)=-\xi_{1}^{4}-2\xi_{1}\xi_{2}^{2}=-\xi_{1}\left(\xi_{1}^{3}+2\xi_{2}^{2}\right).\]
Outside from $\xi_{1}=0$, this determinant vanishes on the curve
give by the equation $\xi_{1}^{3}+2\xi_{2}^{2}=0$. This means that
Condition {}``Iso-energetic'' is not satisfied on this curve.
\end{example*}
\noindent 

\begin{example*}
\noindent \textbf{{}``Rüssmann'' but not {}``Weak''.} On $\mathcal{B}=\mathbb{R}^{2}\setminus0$,
let us consider the function $F\left(\xi\right)=\xi_{1}^{4}+\xi_{2}^{4}$.
The differential is $dF=4\xi_{1}^{3}d\xi_{1}+4\xi_{2}^{3}d\xi_{2}$,
implying that for each $X\in\mathbb{R}^{2}$, one has $\Omega_{X}=dF\left(X\right)=4\xi_{1}^{3}X_{1}+4\xi_{2}^{3}X_{2}$.
The associated resonant set $\Sigma_{X}=\left\{ \left(\xi_{1},\xi_{2}\right)\mid\xi_{1}^{3}X_{1}+\xi_{2}^{3}X_{2}=0\right\} $
is simply a line passing through the origine and with slope equal
to $-\left(\frac{X_{1}}{X_{2}}\right)^{\frac{1}{3}}$, without the
origin point. This set has an empty interior and therefore Condition
{}``Rüssmann'' is satisfied. On the other hand, the differential
of $\Omega_{X}$ is given by $d\Omega_{X}=12\left(\xi_{1}^{2}X_{1}d\xi_{1}+\xi_{2}^{2}X_{2}d\xi_{2}\right)$.
We can see that for some vectors $X$, the differential $d\Omega_{X}$
vanishes at some points belonging to the $\Sigma_{X}$, and thus Condition
{}``Weak'' is not satisfied. For example, for $X=\left(X_{1},0\right)$
the resonant set $\Sigma_{X}$ is the vertical axis $\left\{ \left(0,\xi_{2}\right);\xi_{2}\neq0\right\} $
without the origin point. Now, at each point of this set, one has
$d\Omega_{X}=0$.
\end{example*}

\section{Somes properties of ND hamiltonians}

Let $\left(H,\mathcal{M}\overset{\pi}{\rightarrow}\mathcal{B}\right)$
be a regular CI system composed of a fibration in lagrangian tori
$\mathcal{M}\overset{\pi}{\rightarrow}\mathcal{B}$ together with
a Hamiltonian $H\in C^{\infty}\left(\mathcal{M}\right)$ constant
along the fibers, i.e. $H=F\circ\pi$ for some $F\in C^{\infty}\left(\mathcal{B}\right)$.
As explained in Section \ref{sec_Geometric_setting}, such a fibration
implies the existence of a natural torsion-free and flat connection
on the base space $\mathcal{B}$, the \emph{Duistermaat connection}.
Moreover, one can easily show that the Hamiltonian vector field $X_{H}$
associated with $H=F\circ\pi$ is tangent to the fibration and its
restriction to each torus $\mathcal{M}_{b}$ is an element of $\mathcal{V}_{\nabla}\left(\mathcal{M}_{b}\right)$,
i.e. is parallel with respect to Weintein's connection. Therefore,
on each torus, $X_{H}$ generates a linear dynamics that can be periodic,
resonant or non-resonant. 

The resonance properties are actually well-parametrized by using the
\emph{Resonance bundle} $E$. Indeed, for any given $k\in\Gamma\left(E\right)$,
the previously defined set $\Sigma_{k}=\Omega_{k}^{-1}\left(0\right)$,
with $\Omega_{k}=dF\left(k\right)$, is actually the set of tori on
which the Hamiltonian vector field $X_{H}$ satisfies at least one
resonance relation which reads $\sum_{j=1}^{d}k_{i}X^{j}=0$ in action-angle
coordinates. On such tori, the dynamics is actually confined in $\left(d-1\right)$-dimensionnal
subtori. If the dynamics of $X_{H}$ is periodic on a torus $\mathcal{M}_{b}$,
this means that it satisfies $d-1$ resonance relations, i.e. $b$
belongs to the intersection of resonant sets $\Sigma_{k_{1}},...,\Sigma_{k_{d-1}}$,
where the $k_{j}$ are linearly independent sections of $\Gamma\left(E\right)$.
On the other hand, if $X_{H}$ is ergodic on some torus $\mathcal{M}_{b}$
, this means that $b$ does not belong to any resonant set. 

\begin{lem}
\label{lem_ND_ergodic_dense}If $F\in C^{\infty}\left(\mathcal{B}\right)$
is {}``Rüssmann'' non-degenerate, then the set of tori on which
the dynamics is ergodic is dense in $\mathcal{B}$.
\end{lem}
\begin{proof}
If $\mathcal{M}_{b}$ is an ergodic torus, this means that $b$ does
not belong to any resonant set $\Sigma_{k}$, i.e. $b\in\mathcal{B}\setminus\bigcup_{k}\Sigma_{k}$.
We will show that $\mathcal{B}\setminus\bigcup_{k}\Sigma_{k}$ is
dense in $\mathcal{B}$ by showing that the interior of $\bigcup_{k}\Sigma_{k}$
is empty. Indeed, whenever the function $F$ satisfies Condition {}``Rüssmann'',
then for each non-vanishing $k\in\Gamma\left(E\right)$, the subset
$\Sigma_{k}$ has empty interior. Moreover, $\Sigma_{k}$ is a closed
subset since it is the inverse image of the point $0\in\mathbb{R}$
by the continuous map $dF\left(k\right):\mathcal{B}\rightarrow\mathbb{R}$.
We can then apply the Baire's Theorem (see for e.g. \cite{dixmier})
which insures that $\bigcup_{k}\Sigma_{k}$ has empty interior.
\end{proof}
\begin{lem}
\label{lem_app_ND_2}If the Hamiltonian $H=F\circ\pi$ is {}``Rüssmann''
non-degenerate, then the space of functions $C^{\infty}\left(\mathcal{M}\right)$
constant along the fibers equals to the space of functions $C^{\infty}\left(\mathcal{M}\right)$
which Poisson-commute with $H$. 
\end{lem}
\begin{proof}
Indeed, if a function $A\in C^{\infty}\left(\mathcal{M}\right)$ satisfies
$\left\{ H,A\right\} =0$, then one has $X_{H}\left(A\right)=0$ and
therefore $A$ is constant along the trajectories of $X_{H}$. For
each torus $\mathcal{M}_{b}$ on which the dynamics is ergodic, this
implies that $A$ is constant over this torus. Now, Lemma \ref{lem_ND_ergodic_dense}
insures that when $H$ satisfies Condition {}``Rüssmann'', then
the set of ergodic tori is dense in $\mathcal{B}$. By continuity,
this shows that $F$ is constant along all the fibers $\mathcal{M}_{b}$.
Conversely, if $F$ is constant along all the fibers, then $\left\{ H,F\right\} =X_{H}\left(F\right)=0$
since $X_{H}$ is tangent to the fibration.
\end{proof}
\begin{cor}
\label{coro_app_ND_3}Let $H=F\circ\pi$ be a {}``Rüssmann'' non-degenerate
Hamiltonian and $A,B\in C^{\infty}\left(\mathcal{M}\right)$ two functions.
Then, the following holds \[
\left\{ A,H\right\} =\left\{ B,H\right\} =0\Rightarrow\left\{ A,B\right\} =0.\]

\end{cor}
\begin{proof}
Indeed, if $A$ and $B$ commute with $H$, Lemma \ref{lem_app_ND_2}
implies that $A$ and $B$ are constant along the fibers. Moreover,
since the fibers are lagrangian, $X_{A}$ is tangent to them and therefore
$\left\{ A,B\right\} =X_{A}\left(B\right)=0$.
\end{proof}
\begin{thm}
\label{prop_app_ND_4}If $\left(H,\mathcal{M}\overset{\pi}{\rightarrow}\mathcal{B}\right)$
is a {}``Rüssmann'' non-degenerate C.I system, then $\mathcal{M}\overset{\pi}{\rightarrow}\mathcal{B}$
is the unique fibration such that $H$ is constant along the fibers.
\end{thm}
\begin{proof}
Indeed, suppose $\mathcal{M}\overset{\pi'}{\rightarrow}\mathcal{B}^{'}$
is another fibration such that $H$ is constant along the fibers.
Let $a_{1},...,a_{d}\in C^{\infty}\left(\mathcal{B}^{'}\right)$ be
a local coordinate system in an open subset $\mathcal{O}^{'}\subset\mathcal{B}^{'}$
and $A_{j}=a_{j}\circ\pi^{'}$, $j=1..d$, the pull-back functions.
Since the differentials $dA_{j}$ are linearly independent in $\pi^{-1}\left(\mathcal{O}^{'}\right)$,
the fibers $\left(\pi^{'}\right)^{-1}$ are given by the level-sets
of the $A_{j}$'s. Moreover, we have $\left\{ A_{j},H\right\} =0$
since by hypothesis $H$ is constant along the lagrangian fibration
$\pi^{'}$. Now, Lemma \ref{lem_app_ND_2} implies that the functions
$A_{j}$ must be constant along the fibers of the first fibration
$\pi$, since $H$ is non-degenerate. This means that the fibers $\pi^{-1}$
are included in the level sets of the $A_{j}$'s, and thus are included
in the fibers of $\pi^{'}$. Since the fibers of both fibrations have
the same dimension, they must coincide.
\end{proof}
\bibliographystyle{plain}
\bibliography{/home/roy/math/biblio/biblio_nico}

\end{document}